\def\ifudf#1{\expandafter\ifx\csname #1\endcsname\relax}
\newif\ifpdf \ifudf{pdfoutput}\pdffalse\else\pdftrue\fi
\ifpdf\pdfpagewidth=210mm \pdfpageheight=297mm
 \else\fi
\hsize=159mm \vsize=245mm \hoffset=0pt \voffset=0pt
\parskip=1ex minus .3ex \parindent=0pt
\everydisplay={\textstyle}
\font\tbb=bbmsl10 \font\sbb=bbmsl10 scaled 700
 \font\fbb=bbmsl10 scaled 500
\newfam\bbfam \textfont\bbfam=\tbb
 \scriptfont\bbfam=\sbb \scriptscriptfont\bbfam=\fbb
\def\bb{\fam\bbfam}%
\font\tfk=eufm10 \font\sfk=eufm7 \font\ffk=eufm5
\newfam\fkfam \textfont\fkfam=\tfk
 \scriptfont\fkfam=\sfk \scriptscriptfont\fkfam=\ffk
\def\fk{\fam\fkfam}%
\def\fn[#1]{\font\TmpFnt=#1\relax\TmpFnt\ignorespaces}
\def\em{\expandafter\ifx\the\font\tensl\rm\else\sl\fi}
\def\dt{\number\day.\number\month.\number\year/\the\time}
\def\\{\hfill\break}
\def\R{{\bb R}}
\def\id{{\rm id}}
\def\dts{{D\over ds}^{\!t}}
\def\<#1>{\langle#1\rangle}
\newcount\secNo \secNo=0
\def\section #1\par{\goodbreak\vskip 3ex\noindent
 \global\advance\secNo by 1 \global\eqnNo=0
 {\fn[cmbx10 scaled 1200]#1}\vglue 1ex}
\def\href#1<#2>{\leavevmode
 \ifpdf\pdfstartlink attr {/Border [0 0 0 ]} goto name {#1}\fi
 {#2}\ifpdf\pdfendlink\fi}
\def\label@#1:#2@{\ifudf{#1}
 \expandafter\xdef\csname#1\endcsname{#2}\else
 \errmessage{label #1 already in use!}\fi}
\input\jobname.lbl
\newwrite\lbl \immediate\openout\lbl=\jobname.lbl
\def\@#1:#2@{\ifpdf\pdfdest name {#1} xyz\fi {#2}%
 \immediate\write\lbl{\string\label@#1:#2@}}
\def\:#1:{\href#1<\ifudf{#1}??\else\csname#1\endcsname\fi>}
\newcount\refNo \refNo=0
\def\refitem#1 {\global\advance\refNo by 1
 \item{\@#1:\number\refNo@.}}
\let\pleqno\eqno \newcount\eqnNo \eqnNo=0
\def\eqno#1$${\global\advance\eqnNo by 1
 \pleqno{\rm(\@#1:\number\secNo.\number\eqnNo@)}$$}
\newtoks\title \newtoks\stitle
\newtoks\author \newtoks\sauthor
\newtoks\funding
\title={Semi-discrete isothermic surfaces}
\stitle={Semi-discrete isothermic surfaces}
\author={F Burstall, U Hertrich-Jeromin, C M\"uller, W Rossman}
\sauthor={Burstall et al}
\funding={This work has been partially supported by
 the Austrian Science Fund (FWF) and
 the Japan Society for the Promotion of Science (JSPS)
 through the FWF/JSPS Joint Project grant I1671-N26
 ``Transformations and Singularities''.}
\ifpdf\pdfoutput=1\pdfadjustspacing=1\pdfinfo{%
 /Title (\the\title) /Author (\the\author) /Date (\dt)}\fi
\headline={\ifnum\pageno=1\hfil\else
 {\fn[cmr7]\the\sauthor\hfil\the\stitle}\fi}
\footline={{\fn[cmr7]\hfil-\folio-\hfil}}
\centerline{{\fn[cmbx10 scaled 1440]\the\title}}\vglue .2ex
\centerline{{\fn[cmr7]\the\author}}\vglue 3em plus 3ex
\centerline{\vtop{\hsize=.8\hsize{\bf Abstract.}\enspace
 A Darboux transformation for polarized space curves is
 introduced and its properties are studied, in particular,
 Bianchi permutability.
 Semi-discrete isothermic surfaces are described as sequences
 of Darboux transforms of polarized curves in the conformal
 $n$-sphere and their transformation theory is studied.
 Semi-discrete surfaces of constant mean curvature are studied
 as an application of the transformation theory.
}}\vglue 2em
\centerline{\vtop{\hsize=.8\hsize{\bf MSC 2010.}\enspace
{\it 53A10\/}, {\it 53C42\/}, 53A30, 37K25, 37K35
}}\vglue 1em
\centerline{\vtop{\hsize=.8\hsize{\bf Keywords.}\enspace
isothermic surface; discrete isothermic net;
Calapso transformation; Darboux transformation;
Bianchi permutability; Lawson correspondence;
B\"acklund transformation; constant mean curvature.
}}\vglue 3em

\section 1. Introduction

Integrable discretizations of surfaces or submanifolds are
intimately related to the transformations of the smooth theory
---
permutability theorems govern consistency of the discretization
in this setting.
This line of thought has probably been most clearly formulated
in [\:bosu07a:] or [\:bosu08:]:
``{\em In this setting, discrete surfaces appear as
  two-dimensional layers of multidimensional discrete nets,
  and their transformations correspond to shifts in the
  transversal lattice directions.
  A characteristic feature of the theory is that all lattice
  directions are on equal footing with respect to the defining
  geometric properties.\/}''

For example, given two Darboux transforms of an isothermic
surface, a fourth surface can be constructed by algebraic means
so that the four surfaces form a quadrilateral of surfaces with
Darboux transforms along the edges of the quadrilateral;
moreover, the (spectral) parameters of the transformations on
opposite edges are equal.
Thus, repeated transformations and application of this
permutability theorem generate a lattice of isothermic
surfaces --- trailing the effect on a single point then
generates a $2$-dimensional grid in space:
a discrete surface with similar properties as a smooth
isothermic surface.
In particular, the discrete surfaces obtained in this way
admit very similar transformations to their smooth analogues,
hence generating multi-dimensional lattices that display the
same properties in all lattice directions,
cf [\:bosu07a:]:
``{\em Discrete surfaces and their transformations become
  indistinguishable.\/}''

Our principal aim here is to explore this relation between
analogous smooth and discrete theories in a setting where
the interplay between them becomes most tangible,
that is, we unify the smooth and discrete aspects in a single
geometric object: a semi-discrete isothermic surface.

Thus, to construct the semi-discrete isothermic surfaces
considered in this text, we follow the idea outlined above,
see also [\:caiv98:]:
we obtain a notion of Darboux transformation for space curves
by observing how the Darboux transformation of an isothermic
surface acts on a single curvature line, cf [\:ho09:]
---
then we construct semi-discrete isothermic surfaces as
sequences of Darboux transforms of curves.

The first part of the text is concerned with the analysis
of the Darboux transformation of curves in Euclidean and
conformal $n$-space:
we introduce the Darboux transformation as a special type
of Ribaucour transformation,
cf [\:boje01:] and [\:imdg:, Sect 8.2],
that preserves a quadratic differential as an additional
structure on a curve --- this Darboux transformation acts
on curves equipped with a quadratic differential,
reminiscent of the fact that the Darboux transformation
of an isothermic surface preserves a holomorphic quadratic
differential that is associated to every isothermic surface.
This notion of Darboux transformation for {\em polarized\/}
curves generalizes the notions of [\:ho09:, Sect 2.6].

In particular, Darboux transformations are given by solutions
of a Riccati equation that includes a (spectral) parameter
---
we discuss a linearization of this Riccati equation using
M\"obius geometric techniques, by means of a connection
on a suitable vector bundle, reflecting the conformal
invariance of the transformation.
Darboux transforms are then obtained from parallel sections
of this connection, that is, as solutions of Darboux's linear
system.
A version of the aforementioned Bianchi permutability theorem
for the Darboux transformation of curves is then derived using
standard gauge theoretic arguments.

Semi-discrete isothermic surfaces are introduced in the second
part of the text, as sequences of Darboux transforms of curves.
As detailed above in the fully discrete case, the Darboux
transformation for semi-discrete isothermic surfaces can
then be consistently defined because of the corresponding
Bianchi permutability theorem for the Darboux transformation
of curves.
Further, we discuss how the Christoffel and Calapso
transformations of semi-discrete isothermic surfaces
occur from corresponding transformations for curves and
suitable permutability theorems.

In summary, we do not only obtain the full transformation
theory for semi-discrete isothermic surfaces, but, relying
on the aforementioned linearization of the Riccati equation
that describes the Darboux transformation, we also obtain a
gauge theoretic characterization of semi-discrete isothermic
surfaces, as are known in the smooth and fully discrete cases,
cf [\:bdpp11:].

As an example of the developed transformation theory and gauge
theoretic approach for semi-discrete isothermic surfaces,
we discuss semi-discrete surfaces of constant mean curvature
in space forms, thereby extending results from
[\:muwa13:] and [\:mu15:].
Here we rely on the intimate relation between the Christoffel
transformation and a mixed area that is used to define the
mean curvature of a surface with Gauss map in a space form,
cf [\:bjl14:] and [\:bjr14:] for the discrete case;
a characterization in terms of linear conserved quantities,
cf [\:busa12:] and [\:bjrs08:],
is then a direct consequence of the gauge theoretic
characterization of semi-discrete isothermic surfaces
and its relation with the Christoffel transformation.

{\it Acknowledgements.\/}
We would like to thank
 M Yasumoto
for fruitful and enjoyable discussions around the subject.
\the\funding

\section 2. Darboux transformations of curves

Our idea is that semi-discrete isothermic surfaces should be
obtained by restricting a sequence of Darboux transforms of
an isothermic surface to a fixed curvature line.
We therefore seek a notion of ``Darboux transformation'' for
(decorated) curves in space and shall find it in much the same
way that a ``B\"acklund transformation'' for curves of constant
torsion was found by restricting B\"acklund transformations of
pseudo-spherical surfaces to a single asymptotic line
in [\:caiv98:, Sect 1.1].
In particular, the Darboux transformation will occur as a
special type of Ribaucour transformation,
cf [\:imdg:, Sect 8.2] or [\:boje01:, Sect 3]:
the Ribaucour transformation of hypersurfaces or submanifolds
descends to their curvature lines, where only the enveloping
condition remains a non-trivial condition.

\proclaim\@def.ribaucour:Def 2.1@.
Two curves $x,\hat x:I\to\R^n$ will be said to form
a {\em Ribaucour pair\/} if they envelop a circle congruence,
i.e.,
if tangents at corresponding points $x(s)$ and $\hat x(s)$
are tangent to a common circle $c(s)$.
Either curve of a Ribaucour pair is a
{\em Ribaucour transform\/} of the other curve.

Using a semi-discrete version of the Clifford algebra
cross ratio of [\:ci97:], cf [\:bu06:, Sect 2.4],
we obtain an algebraic characterization of Ribaucour pairs
of curves:

\proclaim\@thm.ribaucour:Def \& Lemma 2.2@.
A pair of curves $(x,\hat x):I\to\R^n\times\R^n$ is
a Ribaucour pair if and only if its {\em tangent cross ratio\/}
$$
  cr := x'\,(x-\hat x)^{-1}\hat x'\,(x-\hat x)^{-1}:I\to\R
\eqno icr$$
is real, where Clifford multiplication is used
for products and inverses.
In that case, the tangent cross ratio is symmetric,
$$
    x'\,(x-\hat x)^{-1}\hat x'\,(x-\hat x)^{-1}
  = \hat x'\,(\hat x-x)^{-1}x'\,(\hat x-x)^{-1}.
$$

Namely, $(x,\hat x)$ envelop a circle congruence if and only
if reflection in the bisector of the secants exchanges
the tangents, that is, if and only if
$$
  x' \parallel (\hat x-x)^{-1}\hat x'(\hat x-x)
  = x'\,{(\hat x-x)^2cr\over x'^2}.
$$
Note that corresponding tangents of a Ribaucour pair of curves
generically intersect in a point $y$, which provides an
alternative method to determine the tangent cross ratio
of a Ribaucour pair from a suitable ratio of scaling factors:
$$
  \hat x + {1\over\hat r}\hat x' = y = x + {1\over r}x',
   \enspace{\rm where}\enspace
  {\hat r\over r} = -{(\hat x-x)^2cr\over x'^2}
\eqno icrpoint$$
since
$$
  \hat x' + {(\hat x-x)^2cr\over x'^2}\,x'
  = \hat x' + (\hat x-x)^{-1}\hat x'(\hat x-x)
  \parallel (\hat x-x).
$$

On the other hand, given a tangent cross ratio function
$cr:I\to\R$, a corresponding Ribaucour transform $\hat x$
of a curve $x:I\to\R^n$ can be recovered as a solution of
the Riccati equation
$$
  \hat x' = cr\,(\hat x-x)(x')^{-1}(\hat x-x).
\eqno ribricc$$
Our next mission will be to linearize this equation by means
of standard M\"obius geometric methods,
see for example [\:imdg:, Chap 1]:
we will replace the Euclidean $\R^n$ as an ambient space by the
conformal $n$-sphere $S^n$, thought of as the projective light
cone of the $(n+2)$-dimensional Minkowski space $\R^{n+1,1}$,
$$
  S^n\cong{\bb P}({\cal L}^{n+1}),
   \enspace{\rm where}\enspace
  {\cal L}^{n+1}=\{y\in\R^{n+1,1}\,|\,(y,y)=0\}.
$$
Thus a (smooth) choice of homogeneous coordinates for a curve
in the conformal $n$-sphere $S^n$ constitutes a {\em lift\/}
of the curve into the light cone ${\cal L}^{n+1}$.
In particular, for lightlike vectors $o,q\in\R^{n+1,1}$ with
inner product $(o,q)=-1$ consider the orthogonal decomposition
$$
  \R^{n+1,1} = \R^{1,1} \oplus_\perp \R^n
  = \<o,q> \oplus_\perp \<o,q>^\perp,
$$
where $\<\dots>$ denotes the span of vectors.
Then the {\em Euclidean lift\/} of an immersed curve
(or, more generally, submanifold)
in $\R^n$ is obtained from the isometry
$$
  \R^n \ni x \mapsto \xi :
  = o + x + {1\over 2}(x,x)\,q \in Q^n :
  = \{y\in{\cal L}^{n+1}\,|\,(y,q)=-1\}
\eqno eucl$$
into the light cone, realizing the Euclidean geometry of $\R^n$
as a {\em subgeometry\/} of the conformal geometry of $S^n$,
cf [\:imdg:, Sect 1.4].
Hence identifying $\R^n\cong Q^n$ as a subset of the conformal
$n$-sphere $S^n$ by means of (\:eucl:), we will also refer to
any light cone map
$$
  \xi:I\to{\cal L}^{n+1}
   \enspace{\rm with}\enspace
  \<\xi>=\<o+x+{1\over 2}(x,x)\,q>
$$
as a {\em lift\/} of a curve $x:I\to\R^n\subset\R^{n+1,1}$
in Euclidean space.
Further,
identifying $\Lambda^2\R^{n+1,1}\cong{\fk o}(\R^{n+1,1})$
via
$$
  \xi\wedge\eta:
  \R^{n+1,1} \ni y \mapsto (\xi\wedge\eta)(y) :
  = (y,\xi)\,\eta - (y,\eta)\,\xi \in \R^{n+1,1},
$$
we linearize (\:ribricc:) by means of a connection
${D\over ds}$ along $x:I\to\R^n\subset\R^{n+1,1}$,
associated to a given tangent cross ratio function
$cr:I\to\R$:

\proclaim\@thm.ribricc:Lemma 2.3@.
A curve $\hat x:I\to\R^n$ is a Ribaucour transform of a curve
$x:I\to\R^n$ with tangent cross ratio function $cr:I\to\R$
if and only if,
for any lifts $\xi$ and $\hat\xi$ of $x$ and $\hat x$,
respectively,
$$
  {D\over ds}\,\hat\xi \parallel \hat\xi,
   \enspace{\sl where}\enspace
  {D\over ds} := {d\over ds}
  - 2cr\,{\xi\wedge\xi'\over(\xi',\xi')}.
\eqno ribsystem$$

Using Euclidean lifts $\xi$ and $\hat\xi$ of two curves
$x$ resp $\hat x$ in $\R^n$, (\:ribsystem:) reads
$$
 \lambda\hat\xi + \hat\xi'
  = {2cr\over(\xi',\xi')}((\hat\xi,\xi)\,\xi'
    - (\hat\xi,\xi')\,\xi)
  = {2cr\,(x-\hat x,x')\over(x',x')}\,\xi
  - {(\hat x-x)^2cr\over x'^2}\,\xi',
$$
where $\lambda={2cr\,(x-\hat x,x')\over(x',x')}$ by taking
inner product with $q$; hence comparing $\R^n$-parts of the
equation we recover (\:icrpoint:), which yields the claim.
Conversely, assuming (\:icrpoint:), expansion of
$$
  0 = (\hat x-x,y-{\hat x+x\over 2})
   \enspace{\rm yields}\enspace
  {D\over ds}\hat\xi = {2cr\,(x-\hat x,x')\over(x',x')}\hat\xi.
$$

Moreover, (\:ribsystem:) is
independent of the choice of lifts $\xi$ resp $\hat\xi$:
multiplication of $\xi$ by a non-zero function does not change
${D\over ds}$ since $(\xi,\xi)\equiv 0$ and
multiplication of $\hat\xi$ by a non-zero function does
obviously not affect (\:ribsystem:) as ${D\over ds}$ satisfies
Leibniz' rule.
This reflects the fact that the Ribaucour transformation
is well defined for curves in the conformal
$S^n\cong{\bb P}({\cal L}^{n+1})$:
on a $1$-dimensional domain the only non-trivial condition
for the Ribaucour transformation is the enveloping condition,
that is, the condition that
$$
  \<\xi,\xi',\hat\xi,\hat\xi'>\cong\R^{2,1}\subset\R^{n+1,1}
$$
define a circle congruence, cf [\:imdg:, Sect 6.6];
the structure of the connection ${D\over ds}$ resp of
(\:ribsystem:) then encodes the tangent cross ratio function:
for arbitrary lifts $\xi$ and $\hat\xi$ of the curves of
a Ribaucour pair (\:ribsystem:) reads
$$
  \hat\xi' = 2cr\,{(\xi,\hat\xi)\over(\xi',\xi')}\,\xi'
  \bmod\<\xi,\hat\xi>.
\eqno icrlift$$
As ${D\over ds}$ is a metric connection, every initial point
$\<\hat\xi(s_0)>\in S^n={\bb P}({\cal L}^{n+1})$ yields a
unique Ribaucour transform of a given curve $\<\xi>:I\to S^n$
via a ${D\over ds}$-parallel section
$\hat\xi:I\to{\cal L}^{n+1}$.
A similar construction of Ribaucour transforms via parallel
sections of a connection is possible in higher dimensions,
though further conditions on the tangent cross ratios in
curvature directions and initial point need to be satisfied.

Clearly, the tangent cross ratio of a Ribaucour pair
$(x,\hat x)$ of curves is parameter dependent.
Considering differentials instead of derivatives we obtain
the {\em infinitesimal cross ratio\/}
$$
  x'ds\,(\hat x-x)^{-1}\hat x'ds\,(\hat x-x)^{-1} = cr\,ds^2,
$$
a parameter independent quadratic differential.
This resonates well with the Darboux transformation of
an isothermic surface preserving a holomorphic quadratic
differential (polarization) that occurs from a factorization
of its Hopf differential,
see [\:imdg:, Lemma 5.2.12, \S5.4.14]
and [\:bu06:, Lemma 2.1, Thm 2.7].
Thus to introduce the Darboux transformation for curves we will
consider a (non-vanishing) reference quadratic differential
({\em polarization\/}) ${ds^2\over m}$ as part of the data.

\proclaim\@def.darboux:Def 2.4@.
A Ribaucour pair of curves $x,\hat x:(I,{ds^2\over m})\to\R^n$
on a {\em polarized domain\/} $(I,{ds^2\over m})$,
$m:I\to\R^\times$,
will be called a {\em Darboux pair\/}
if its infinitesimal cross ratio is a constant multiple
of the reference polarization,
$$
  x'(\hat x-x)^{-1}\hat x'(\hat x-x)^{-1} = {\mu\over m}
   \enspace{\rm with}\enspace
  \mu \equiv const\in\R.
\eqno dcrs$$
Either curve of a Darboux pair is a
{\em Darboux transform\/} of the other polarized curve.

Clearly, since the tangent cross ratio of a Ribaucour
pair $(x,\hat x)$ is symmetric in $x$ and $\hat x$,
so is the notion of Darboux transformation.
An example of the Darboux transformation in Euclidean geometry
is given by a tractrix construction, cf [\:ho09:, Def 2.41]:
if $y:I\to\R^n$ denotes an arc-length parametrized curve
then
$$
  x_\pm := y \pm {1\over 2\sqrt\mu}\,y',
   \enspace{\rm where}\enspace
  \mu\in(0,\infty),
\eqno tractrix$$
form a Darboux pair of curves with respect to their common
{\em arc-length polarization\/} ${ds^2\over m}=|dx_\pm|^2$;
namely, with a principal normal field $n:I\to S^{n-1}$ of $y$
and the corresponding curvature $\kappa$ we obtain
$$\left.\matrix{
   x_\pm' = y' \pm {\kappa\over 2\sqrt\mu}\,n \hfill\cr
   x_+-x_- = {1\over\sqrt\mu}\,y' \hfill\cr
  }\right\}\enspace\Rightarrow\enspace
  cr\,ds^2 = \mu\,|x_\pm'ds|^2.
$$

Given a polarized curve $x:(I,{ds^2\over m})\to\R^n$ every
Darboux transform $\hat x:(I,{ds^2\over m})\to\R^n$ is
obtained as a solution of a Riccati equation
$$
  \hat x' = \mu\,(\hat x-x)(mx')^{-1}(\hat x-x);
\eqno riccati$$
as for the Ribaucour transformation the Darboux transformation
is most naturally considered as a transformation of polarized
curves in the conformal $n$-sphere $S^n={\bb P}({\cal L}^{n+1})$,
where the Riccati equation (\:riccati:) linearizes and becomes
{\em Darboux's linear system\/},
cf [\:da99:],
that is, yields the zero curvature representation of the
associated integrable system.
In particular, for a polarized curve
$$
  \<\xi>:(I,{ds^2\over m})\to S^n={\bb P}({\cal L}^{n+1}),
   \enspace{\rm where}\enspace
  \xi:I\to\R^{n+1,1}
$$
denotes {\em any\/} light cone lift of the curve, we introduce
a $1$-parameter family of (flat) connections by
$$
  \dts := {d\over ds}
  - {2t\over m}\,{\xi\wedge\xi'\over(\xi',\xi')}.
\eqno dtsmooth$$
Note that, as $\xi$ takes values in the light cone,
the connections $\dts$ do not depend on the choice of lift.
Now \:thm.ribricc: for the Ribaucour transformation directly
yields a similar assertion for the Darboux transformation:

\proclaim\@thm.darboux:Def \& Cor 2.5@.
The family $(\dts)_{t\in\R}$ from {\rm(\:dtsmooth:)} will be
referred to as the {\em isothermic family of connections\/}
of a polarized curve.  The Darboux transforms of a polarized
curve $\<\xi>:(I,{ds^2\over m})\to S^n$, with respect to a
parameter $\mu\in\R$, are given by $\dts|_{t=\mu}$-parallel
sections, that is,
by solutions of {\em Darboux's linear system\/}
$$
  \dts\big|_{t=\mu}\hat\xi = 0.
\eqno dls$$

Consequently, any polarized curve admits a $(1+n)$-parameter
family of Darboux transforms as any choice of the spectral
parameter $t=\mu$ and of an initial point in $S^n$ yields
a unique Darboux transform via a $\dts|_{t=\mu}$-parallel
light cone section.
On the other hand, this new description (\:dls:) of the
Darboux transformation lacks the symmetry of (\:dcrs:):
thus our next goal is to understand how Darboux's linear
system changes under the Darboux transformation.

\proclaim\@thm.dtsgauge:Lemma 2.6@.
If $\<\xi>,\<\hat\xi>:(I,{ds^2\over m})\to S^n$ form
a Darboux pair, then their isothermic loops of connections
are related by a gauge transformation
$$
  \hat\dts = \Gamma_{\<\xi>}^{\<\hat\xi>}(1-{t\over\mu})
   \cdot\dts, \enspace{\sl where}\enspace
  \Gamma_{\<\xi>}^{\<\hat\xi>}(r)y := \cases{
   ry & for $y\in\<\hat\xi>$, \cr
   y & for $y\in\<\xi,\hat\xi>^\perp$, \cr
   y/r & for $y\in\<\xi>$. \cr}
$$

Note that $\Gamma_{\<\xi>}^{\<\hat\xi>}(r)\in O(\R^{n+1,1})$
realizes the M\"obius geometric cross ratio of four points
on a circle:
$$
  cr(\<\hat\xi>,\<\eta>,\<\xi>,\<\zeta>) = r
   \enspace\Leftrightarrow\enspace
  \<\zeta> = \<\Gamma_{\<\xi>}^{\<\hat\xi>}(r)\,\eta>;
\eqno dcr$$
namely, fixing relative normalizations of $\xi$ and $\hat\xi$
so that $(\xi,\hat\xi)=-1$ and using $\xi\simeq o$ and
$\hat\xi\simeq q=\infty$ as origin and point at infinity
for Euclidean lifts (\:eucl:) of $\<\eta>$ and $\<\zeta>$,
$$
  \eta = \xi + y + {1\over 2}(y,y)\,\hat\xi
   \enspace{\rm and}\enspace
  \zeta = \xi + z + {1\over 2}(z,z)\,\hat\xi,
$$
we find, by computation of the (cyclic) Clifford algebra
cross ratio in $\R^n$ of [\:ci97:], that
$$
  cr(\<\hat\xi>,\<\eta>,\<\xi>,\<\zeta>) = y^{-1}z = r
   \enspace\Leftrightarrow\enspace
  \zeta = r\,\Gamma_{\<\xi>}^{\<\hat\xi>}(r)\,\eta.
$$
Changing viewpoint, (\:dcr:) yields a parametrization
of the circumcircle of three points,
$\<\xi>,\<\hat\xi>,\<\eta>\in S^n$,
by cross ratio:
$$
  \R\cup\{\infty\}\ni r\mapsto
  \<\Gamma_{\<\xi>}^{\<\hat\xi>}(r)\,\eta>\in S^n.
$$

\:thm.dtsgauge: follows directly from the following,
more general statement,
cf [\:bdpp11:, Thm 3.10 and Prop 3.11]:

\proclaim\@thm.gauge:Lemma 2.7@.
Let $\<\xi>,\<\hat\xi>:I\to S^n$ be complementary,
$\<\xi(s)>\neq\<\hat\xi(s)>$ for all $s\in I$, and
$\gamma\in\Omega^1(\<\xi>\wedge\<\xi>^\perp)$ and
$\hat\gamma\in\Omega^1(\<\hat\xi>\wedge\<\hat\xi>^\perp)$
$1$-forms with values in $\<\xi>\wedge\<\xi>^\perp$ resp
$\<\hat\xi>\wedge\<\hat\xi>^\perp\subset{\fk o}(\R^{n+1,1})$.
Then
$$
  (d+t\hat\gamma)=\Gamma_{\<\xi>}^{\<\hat\xi>}(1-{t\over\mu})
   \cdot(d+t\gamma) \enspace\Leftrightarrow\enspace \cases{
  (d+\mu\gamma)\,\hat\xi\parallel\hat\xi & and \cr
  (d+\mu\hat\gamma)\,\xi\parallel\xi. \cr}
$$

We decompose $\R^{n+1,1}
 = \<\xi>\oplus\<\hat\xi>\oplus\<\xi,\hat\xi>^\perp$
and use the corresponding projections
$$
  \pi:\R^{n+1,1}\to\<\xi>, \enspace
  \hat\pi:\R^{n+1,1}\to\<\hat\xi>
   \enspace{\rm and}\enspace
  \varpi:\R^{n+1,1}\to\<\xi,\hat\xi>^\perp,
$$
so that
$$
  \id = \pi + \varpi + \hat\pi
   \enspace{\rm and}\enspace
  \Gamma_{\<\xi>}^{\<\hat\xi>}(r)
  = {1\over r}\pi + \varpi + r\hat\pi;
$$
note that, since $\xi$ and $\hat\xi$ are lightlike,
for example $\pi(v)={(v,\hat\xi)\over(\xi,\hat\xi)}\,\xi$
for $v\in\R^{n+1,1}$.
Thus we have
$$
  \gamma = \pi\circ\gamma\circ\varpi + \varpi\circ\gamma\circ\hat\pi
   \enspace{\rm and}\enspace
  \hat\gamma=\hat\pi\circ\gamma\circ\varpi+\varpi\circ\gamma\circ\pi
$$
as $\gamma\in\Omega^1(\<\xi>\wedge\<\xi>^\perp)$ yields:
{\parindent=2em
\item{$\bullet$} $\gamma\circ\pi=0$
 since $\gamma(\xi)=0$;
\item{$\bullet$} $\gamma\circ\varpi=\pi\circ\gamma\circ\varpi$
 since $\gamma(\<\xi>^\perp)\subset\<\xi>$;
\item{$\bullet$} $\gamma\circ\hat\pi
  =(\pi+\varpi)\circ\gamma\circ\hat\pi$
 since $\gamma(\R^{n+1,1})\subset\<\xi>^\perp$;
 further
\item{$\bullet$} $\pi\circ\gamma\circ\hat\pi=0$
 since $\gamma\in\Omega^1({\fk o}(\R^{n+1,1}))$ is skew symmetric,
 so that $\gamma(\hat\xi)\perp\hat\xi$.
\par}
Similarly the derivative decomposes as
$$
  d = D - \beta - \hat\beta,
   \enspace{\rm where}\enspace\cases{
  \beta :
   = -(\hat\pi\circ d\circ\varpi+\varpi\circ d\circ\pi)
   \in\Omega^1(\<\hat\xi>\wedge\<\xi,\hat\xi>^\perp) &
   and \cr
  \hat\beta :
   = -(\pi\circ d\circ\varpi+\varpi\circ d\circ\hat\pi)
   \in\Omega^1(\<\xi>\wedge\<\xi,\hat\xi>^\perp) \cr}
$$
yield the derivatives of the maps
$\<\xi>,\<\hat\xi>:I\to S^n\subset{\bb P}(\R^{n+1,1})$,
as $(\xi,d\xi)=0$ implies $\hat\pi\circ d\circ\pi=0$
so that $d\xi=-\beta\xi\bmod\xi$ and similarly for $\hat\xi$,
and
$$
  D = \pi\circ d\circ\pi
    + \hat\pi\circ d\circ\hat\pi
    + \varpi\circ d\circ\varpi
$$
is a metric connection on the vector bundle
$\<\xi>\times\<\hat\xi>\times\<\xi,\hat\xi>^\perp$.
Hence we compute
$$\matrix{
  \Gamma_{\<\xi>}^{\<\hat\xi>}(r)\cdot(d+t\gamma)
   -  (d+t\hat\gamma)
  &=& (D - r\beta + {1\over r}(t\gamma-\hat\beta))
   -  (D + (t\hat\gamma-\beta) - \hat\beta) \hfill\cr
  &=& {1\over r}(t\gamma-(1-r)\hat\beta)
   -  (t\hat\gamma-(1-r)\beta). \hfill\cr}
$$
Consequently, with $r=1-{t\over\mu}$, we obtain the claimed
result:
$$
  \Gamma_{\<\xi>}^{\<\hat\xi>}(1-{t\over\mu})\cdot(d+t\gamma)
  = (d+t\hat\gamma)
  \enspace\Leftrightarrow\enspace\left\{\matrix{
   \mu\gamma = \hat\beta \cr
   \mu\hat\gamma = \beta \cr
  }\right\}\enspace\Leftrightarrow\enspace\cases{
   (D-\beta+(\mu\gamma-\hat\beta))\,\hat\xi=0 \bmod\hat\xi, \cr
   (D+(\mu\hat\gamma-\beta)-\hat\beta)\,\xi=0 \bmod\xi. \cr}
$$

Note that we imposed virtually no conditions on the maps
$\<\xi>$ and $\<\hat\xi>$, in particular, we used no
assumptions on the dimension of their common domain,
nor on the $1$-forms $\gamma$ and $\hat\gamma$, apart from
the fact that they take values in the same places as
the derivatives of $\<\xi>$ resp $\<\hat\xi>$ do.
This reflects the utmost generality of the statement,
cf [\:bdpp11:].

It is now readily derived that the well known Bianchi
permutability theorem for the Darboux transformation of
isothermic surfaces descends to a permutability theorem for
the Darboux transformation of polarized curves,
cf [\:bi04:, \S3] or [\:bdpp11:, Sect 4.2]:

\proclaim\@thm.biquad:Thm 2.8@.
Given Darboux transforms
$\<\xi_0>,\<\xi_1>:(I,{ds^2\over m})\to S^n$
of a polarized curve $\<\xi>:(I,{ds^2\over m})\to S^n$
for different spectral parameters $\mu_0\neq\mu_1$,
there is a simultaneous Darboux transform $\<\xi_{01}>$
of $\<\xi_i>$ $(i=0,1)$ with parameters $\mu_{1-i}$,
respectively,
given algebraically by a constant cross ratio function
$$
  cr(\<\xi>,\<\xi_0>,\<\xi_{01}>,\<\xi_1>)
  \equiv {\mu_1\over\mu_0}.
$$

By symmetry it is sufficient to show that $\<\xi_{01}>$
is a $\mu_1$-Darboux transform of $\<\xi_0>$,
using a cross ratio identity, cf [\:imdg:, \S4.9.11].
Thus we determine $\<\xi_{01}>$ from
$$
  cr(\<\xi_0>,\<\xi_1>,\<\xi>,\<\xi_{01}>)
  = 1 - {\mu_1\over\mu_0},
   \enspace{\rm hence}\enspace
  \<\xi_{01}>
  = \<\Gamma^{\<\xi_0>}_{\<\xi>}(1-{\mu_1\over\mu_0})\,\xi_1>
$$
by (\:dcr:);
assuming that $\xi_0$ is a $\dts|_{t=\mu_0}$-parallel section,
$\dts|_{t=\mu_0}\xi_0=0$, and using \:thm.dtsgauge: we obtain
$$
  {(\dts)}_0\big|_{t=\mu_1}\xi_{01}
  = (\Gamma^{\<\xi_0>}_{\<\xi>}(1-{\mu_1\over\mu_0})
    \cdot\dts\big|_{t=\mu_1})
    (\Gamma^{\<\xi_0>}_{\<\xi>}(1-{\mu_1\over\mu_0})\,\xi_1)
  = \Gamma^{\<\xi_0>}_{\<\xi>}(1-{\mu_1\over\mu_0})
    (\dts\big|_{t=\mu_1}\,\xi_1)
  = 0.
$$

The relation between the isothermic loops of connections
in a ``{\em Bianchi quadrilateral\/}'' is now obtained
from \:thm.dtsgauge:;
in fact, a more abstract statement holds true, that yields the
key to an analogous approach to discrete isothermic nets,
cf [\:bjrs08:, Lemma 2.6] or [\:bdpp11:, Lemma 4.7]:
$$
    \Gamma^{\<\xi_{01}>}_{\<\xi_0>}(1-{t\over\mu_1})
    \circ\Gamma^{\<\xi_0>}_{\<\xi>}(1-{t\over\mu_0})
  = \Gamma^{\<\xi_1>}_{\<\xi_0>}({1-t/\mu_1\over 1-t/\mu_0})
  = \Gamma^{\<\xi_{01}>}_{\<\xi_1>}(1-{t\over\mu_0})
    \circ\Gamma^{\<\xi_1>}_{\<\xi>}(1-{t\over\mu_1}).
\eqno bigauge$$

As for isothermic surfaces, cf [\:imdg:, \S5.6.8],
or isothermic submanifolds, cf [\:bdpp11:, Sect 4.3],
the ``{\em Bianchi cube\/}'' permutability theorem for the
Darboux transformation of polarized curves is now obtained
in a completely algebraic way from \:thm.biquad: using
(\:bigauge:) or the hexahedron lemma of [\:imdg:, \S4.9.13],
that is, the ``3D consistency of discrete isothermic nets'',
see [\:bosu08:, Thm 4.26]:

\proclaim\@thm.bicube:Thm 2.9@.
Given three Darboux transforms $\<\xi_i>$, $i=0,1,2$,
of a polarized curve $\<\xi>:(I,{ds^2\over m})\to S^n$
with different parameters, $\mu_i\neq\mu_j$ for $i\neq j$,
there is a simultaneous Darboux transform $\<\xi_{012}>$,
with parameters $\mu_k$,
of the simultaneous Darboux transforms $\<\xi_{ij}>$
of $\<\xi_i>$ and $\<\xi_j>$,
where $i,j,k\in\{0,1,2\}$ are pairwise distinct.

\section 3. Semi-discrete isothermic surfaces

We are now prepared to define semi-discrete isothermic
surfaces, as sequences of Darboux transforms of curves
---
where we will parametrize the sequence by a graph $G=(V,E)$ and
oriented edges $(ij)\in E$ will be denoted by the ordered pair
of their endpoints $i,j\in V$.

\proclaim\@def.sdi:Def 3.1@.
A map $\<\xi>:\Sigma\to S^n$, $(i,s)\mapsto\<\xi_i(s)>$,
on a semi-discrete domain $\Sigma=G\times I$ will be called
a {\em semi-discrete isothermic surface\/} if there is a
polarization ${ds^2\over m}$ on $I$ so that any adjacent
curves $\<\xi_i>$ and $\<\xi_j>$, $(ij)\in E$, form a
Darboux pair of curves on the polarized interval
$(I,{ds^2\over m})$.

Apart from a sign issue this definition recovers the
semi-discrete isothermic surfaces of [\:muwa13:].
With a stereographic projection
$$
  x:\Sigma\to\R^n, \enspace
  \<\xi> = \<o+x+{1\over 2}(x,x)\,q>,
$$
of $\<\xi>$ the condition of adjacent curves forming Darboux
pairs can be formulated using the cross ratio condition
(\:dcrs:),
$$
  x_i'(x_j-x_i)^{-1}\,x_j'(x_j-x_i)^{-1} = {\mu_{ij}\over m}.
$$
Then, reality of the cross ratio is equivalent to two adjacent
curves enveloping a circle congruence while squaring the (real)
cross ratio yields
$$
  \mu_{ij}^2
  = {m^2(x_i',x_i')(x_j',x_j')\over(x_j-x_i)^4}
  = ({\nu_i\nu_j\over(x_j-x_i,x_j-x_i)})^2,
   \enspace{\rm where}\enspace
  \nu:=\sqrt{m(x',x')}.
\eqno sdicr$$
Note that we need to assume $m>0$ here;
$\mu_{ij}<0$ then yields the semi-discrete isothermic surfaces
of [\:muwa13:, Cor 4.2], cf [\:muwa13:, Lemma 3.3]
with $\sigma=-{1\over\mu}$ and $\tau={1\over m}$.
Below we drop this restriction on the sign of the tangent cross
ratio of adjacent curves and also allow Darboux transformations
such as the ones from (\:tractrix:),
or those of the generators of a cylinder to the corresponding
curvature lines of a bubbleton,
cf [\:jepe97:, Sect 9].

As in the smooth and fully discrete cases,
cf [\:bosu08:, Thms 1.32 \& 4.31],
a semi-discrete isothermic surface can be characterized by
the existence of a {\em Moutard lift\/}
$$
  \xi:\Sigma\to{\cal L}^{n+1},
   \enspace{\rm where}\enspace
  A(\xi,\xi)_{ij} := \xi_{ij}'\wedge d_{ij}\xi = 0
   \enspace{\rm with}\enspace\cases{
  d_{ij}\xi := \xi_j-\xi_i, \cr
  \xi_{ij} := {1\over 2}(\xi_i+\xi_j) \cr}
\eqno moutard$$
denoting the discrete derivative resp edge function of $\xi$,
that is, the existence of a light cone lift $\xi$ with vanishing
(algebraic) area element $A(\xi,\xi)ds=0$.
Namely, if $\xi$ satisfies (\:moutard:) then, by (\:icrlift:),
$$
  \xi_{ij}' \parallel d_{ij}\xi = 0 \bmod \<\xi_i,\xi_j>
   \enspace\Rightarrow\enspace
  cr = -{(\xi_i',\xi_i')\over 2(\xi_i,\xi_j)}
     = -{(\xi_j',\xi_j')\over 2(\xi_i,\xi_j)},
$$
in particular, $d_{ij}(\xi',\xi')=0$;
further
$
  (\xi_i,\xi_j)'
  = 2(\xi_{ij},\xi_{ij})'
  = 4(\xi_{ij},\xi_{ij}')
  = 0,
$
as $\xi_{ij}'\parallel d_{ij}\xi$ again,
showing that $\<\xi>$ is isothermic with
$$
  cr = {\mu\over m},
   \enspace{\rm where}\enspace\cases{
  m := {1\over(\xi',\xi')} & so that $d_{ij}m=0$, \cr
  \mu_{ij} := {1\over(d_{ij}\xi,d_{ij}\xi)} &
   so that $\mu_{ij}'=0$. \cr}
$$
Conversely, if $\<x>$ is isothermic with $cr={\mu\over m}$,
where $m>0$ without loss of generality, then we obtain a
Moutard lift by
$$
  \xi := \pm{x\over\sqrt{m(x',x')}},
   \enspace{\rm where}\enspace
  \mu_{ij}(\xi_i,\xi_j)<0
$$
fixes the sign of $\xi$ up to a global sign choice
---
clearly $m(\xi',\xi')\equiv 1$, hence using (\:icrlift:) twice
we learn that
$
  \xi_j'
  = 2\mu(\xi_i,\xi_j)\xi_i'
  = 4\mu^2(\xi_i,\xi_j)^2\xi_j'
  \bmod \<\xi_i,\xi_j>
$
and, in particular, $(\xi_i,\xi_j)=-{1\over2\mu_{ij}}$;
consequently,
$$
  \xi_{ij}' = 0 \bmod \<\xi_i,\xi_j>
   \enspace{\rm and}\enspace
  \xi_{ij}' \perp \xi_{ij}
   \enspace{\rm hence}\enspace
  \xi_{ij}' \parallel d_{ij}\xi.
$$

Our definition has the full isothermic transformation theory
built in: clearly, it directly lends itself to define the
Darboux transformation of semi-discrete isothermic surfaces
via the Darboux transformation of polarized curves and the
Bianchi permutability theorem \:thm.biquad:,
showing that the Darboux transformations of curves
``fit together'' to form a new semi-discrete isothermic
surface.
The Bianchi cube theorem \:thm.bicube: then yields the
semi-discrete analogue of the usual Bianchi permutability
theorem \:thm.biquad:, cf [\:bi04:].

Note that the curves of the original semi-discrete isothermic
surface and of its Darboux transforms become indistinguishable
or, otherwise said,
a Darboux transform leads to a single ``larger'' semi-discrete
isothermic surface with
$$
  \tilde\Sigma = \tilde G\times I, \enspace
  \tilde V = V\mathop{\dot\cup}V \enspace{\rm and}\enspace
  \tilde E = E\mathop{\dot\cup}V\mathop{\dot\cup}E
$$
as a domain, that is, with two copies of the original graph $G$
and ``vertical'' edges added between corresponding points of
the two copies.

Similarly, a ``Christoffel duality'' for polarized curves in
$\R^n$, together with a suitable permutability theorem, yields
a Christoffel duality for semi-discrete isothermic surfaces,
cf [\:muwa13:, Thm 4.3],
as long as there are no cycles in the discrete part $G$ of
the domain $\Sigma$ of the semi-discrete isothermic surface.

\proclaim\@def.chr:Def 3.2@.
Two curves $x,x^\ast:(I,{ds^2\over m})\to\R^n$ on a polarized
domain will be said to be {\em Christoffel dual\/} if
$$
  dx\,dx^\ast = x'(x^\ast)'ds^2 = {ds^2\over m}.
$$

This duality for polarized curves yields the first equation
of [\:muwa13: (10)], the second yields permutability with the
Darboux transformation:

\proclaim\@thm.cdquad:Thm 3.3@.
Let $x:(I,{ds^2\over m})\to\R^n$ be a polarized curve with
Christoffel dual $x^\ast$ and a Darboux transform $\hat x$,
that is, $\hat x'=\mu\,(\hat x-x)(x^\ast)'(\hat x-x)$;
then
$$
  \hat x^\ast := x^\ast + {1\over\mu(\hat x-x)}
$$
is simultaneously Christoffel dual to $\hat x$ and a Darboux
transform of $x^\ast$.

It is straightforward to verify the claim, using that
Christoffel duality is involutive and by computing
$$
  (m\hat x')^{-1}
  = (\hat x^\ast)'
  = \mu(\hat x^\ast-x^\ast)\,x'(\hat x^\ast-x^\ast).
$$
This permutability theorem yields a version of the consistency
check in [\:muwa13:] --- thus it shows that Christoffel duals
of consecutive curves of a semi-discrete isothermic surface can
be positioned to form a new semi-discrete isothermic surface,
its Christoffel dual, as long as the original surface has no
discrete cycles.
In fact, as long as $m>0>\mu$ our Christoffel duality is the
same as the one of [\:muwa13:], as the defining equations
coincide:
$$
  (x^\ast)' = {1\over mx'}
   \enspace{\rm and}\enspace
  d_{ij}x^\ast = {1\over\mu_{ij}d_{ij}x}.
\eqno christoffel$$
As in the smooth and fully discrete cases,
cf [\:imdg:, \S5.3.12] resp [\:bopi96:, Sect 7],
this Christoffel transformation gives rise to a ``Weierstrass
representation'' for semi-discrete minimal surfaces,
cf [\:muwa13:, Sect 5] and [\:roya12:].

To obtain the Calapso transformation,
or ``conformal deformation'',
of a polarized curve we introduce a (family of) gauge
transformation(s) that trivialize the (flat) connections
of its isothermic family of connections (\:dtsmooth:):

\proclaim\@def.calapso:Def 3.4@.
Let $\<\xi>:(I,{ds^2\over m})\to S^n$ be a polarized curve
with isothermic family of connections $\dts$;
$$
  T^t:(I,{ds^2\over m})\to O(\R^{n+1,1})
   \enspace{\sl with}\enspace
  {d\over ds}\circ T^t = T^t\circ\dts
   \enspace{\sl for}\enspace
  t\in\R^\times
\eqno cgt$$
are called the {\em Calapso transformations\/} of $\<\xi>$,
each curve $\<\xi^t>:=\<T^t\xi>$ is a {\em Calapso transform\/}
of $\<\xi>$.

Note that the Calapso transformations $T^t$ can be chosen to
take values in the orthogonal group since the connections
$\dts$ are metric connections;
then they are unique up to post-composition by a (lift of a)
M\"obius transformation.

To carry this Calapso deformation for curves over to one
for semi-discrete isothermic surfaces we will again require
a permutability theorem;
to this end, we will need to gain control of the Calapso
transformations of the Calapso and Darboux transforms of
a polarized curve.

Thus let $\<\tilde\xi>$ denote a Calapso transform of a
polarized curve $\<\xi>$, that is, $\tilde\xi=T^\tau\xi$
for some $\tau\in\R$.
Then (\:cgt:) yields $\tilde\xi'=T^\tau\xi'$, hence
$$
  \tilde\dts\circ T^\tau
  = ({d\over ds}-{2t\over m}
    {\tilde\xi\wedge\tilde\xi'\over(\tilde\xi',\tilde\xi')})
    \circ T^\tau
  = T^\tau\circ({D\over ds}^{\!\tau}-{2t\over m}
    {\xi\wedge\xi'\over(\xi',\xi')})
  = T^\tau\circ{D\over ds}^{\!\tau+t}
$$
so that
$$
  \tilde T^tT^\tau\circ{D\over ds}^{\!\tau+t}
  = \tilde T^t\circ\tilde\dts\circ T^\tau
  = {d\over ds}\circ\tilde T^tT^\tau
   \enspace\Rightarrow\enspace
  \tilde T^tT^\tau = T^{\tau+t}.
$$

For a Darboux transformation $\<\hat\xi>$ of $\<\xi>$, where
$\hat\xi$ satisfies (\:dls:), $\dts\big|_{t=\mu}\hat\xi=0$,
we apply \:thm.dtsgauge: to obtain
$$
  \hat T^t\Gamma_{\<\xi>}^{\<\hat\xi>}(1-{t\over\mu})\circ\dts
  = \hat T^t\circ\hat\dts\circ
    \Gamma_{\<\xi>}^{\<\hat\xi>}(1-{t\over\mu})
  = {d\over ds}\circ\hat T^t
    \Gamma_{\<\xi>}^{\<\hat\xi>}(1-{t\over\mu})
   \enspace\Rightarrow\enspace
  \hat T^t\Gamma_{\<\xi>}^{\<\hat\xi>}(1-{t\over\mu}) = T^t.
$$

Thus we have proved the following:

\proclaim\@thm.ttd:Lemma 3.5@.
Up to post-composition by a M\"obius transformation:
{\parindent=2em
\item{\rm(i)} $\tilde T^tT^\tau=T^{\tau+t}$ for a
 Calapso transform $\<\tilde\xi>=\<T^\tau\xi>$ of $\<\xi>$;
\item{\rm(ii)}
 $\hat T^t\Gamma_{\<\xi>}^{\<\hat\xi>}(1-{t\over\mu})=T^t$
 for a Darboux transform $\<\hat\xi>$ of $\<\xi>$,
 where $\dts\big|_{t=\mu}\hat\xi=0$.
}

In particular, if $\<\hat\xi>$ is a Darboux transform
of $\<\xi>$ with parameter $\mu$,
that is, $\dts\big|_{t=\mu}\hat\xi=0$,
then we learn from \:thm.ttd: (ii) that
$$
  \<\hat T^\tau\hat\xi> = \<T^\tau\hat\xi>
  \enspace{\rm for}\enspace\tau\neq\mu;
$$
hence, by \:thm.ttd: (i),
$$
  \<\tilde T^{\mu-\tau}\hat T^\tau\hat\xi>
  = \<\tilde T^{\mu-\tau}T^\tau\hat\xi>
  = \<T^\mu\hat\xi>
  \equiv const,
$$
showing that $\<T^\tau\hat\xi>$ yields simultaneously
a Calapso transform of $\<\hat\xi>$ and
a Darboux transform of $\<T^\tau\xi>$.
Thus we obtain the desired permutability theorem for the
Darboux and Calapso transformations, which ensures that
the Calapso transformation for curves extends to one
for semi-discrete isothermic surfaces,
as long as the surfaces do not have discrete cycles.

\proclaim\@thm.tdquad:Thm 3.6@.
Suppose that $\<\xi>,\<\hat\xi>:(I,{ds^2\over m})\to S^n$
form a Darboux pair with parameter $\mu\in\R^\times$,
and let $\tau\neq\mu$.
Then their Calapso transforms $\<T^\tau\xi>$ and
$\<\hat T^\tau\hat\xi>$ form, if suitably positioned,
a Darboux pair with parameter $\mu-\tau$.

Thus we have learned that the transformations of curves
extend to sequences of Darboux transforms of curves,
that is, to semi-discrete isothermic surfaces,
by means of suitable permutability theorems:

\proclaim\@thm.sdtrafos:Thm \& Def 3.7@.
The Darboux, Christoffel and Calapso transformations for
polarized curves extend to corresponding transformations
of cycle-free semi-discrete isothermic surfaces.

An alternative approach to the transformations of semi-discrete
isothermic surfaces is more directly based on a family of
semi-discrete flat connections,
similar to those of \:thm.darboux::

\proclaim\@thm.sdiconn:Def \& Thm 3.8@.
A {\em semi-discrete connection\/} on a vector bundle $X$
over a domain $\Sigma=G\times I$, $G=(V,E)$, is a pair
$(\Gamma,\nabla)$, consisting of
connections over each smooth component,
 $\nabla_i$ on $X_i$ for each $i\in V$,
and vector bundle isomorphisms between components,
 $\Gamma_{ij}:X_j\to X_i$ for each $(ij)\in E$;
a semi-discrete connection is {\em flat\/} if all $\nabla_i$
are flat and are gauge equivalent via $\Gamma_{ij}$,
$$
  \forall i\in V:R^{\nabla_i} = 0
   \enspace{\sl and}\enspace
  \forall(ij)\in E:\Gamma_{ij}\cdot\nabla_j = \nabla_i.
$$
A semi-discrete surface $\<\xi>:\Sigma\to S^n$ is isothermic if
and only if there are functions $m$ and $\mu$ on $I$ and $E$,
respectively,
so that the associated {\em isothermic loop of connections\/}
$(\Gamma^t,\nabla^t)_{t\in\R}$ consists of flat connections,
where
$$
  \Gamma^t_{ij} :
  = \Gamma^{\<\xi_i>}_{\<\xi_j>}(1-{t\over\mu_{ij}})
   \enspace{\sl and}\enspace
  \nabla^t_i :
  = \dts\big|_i
  = {d\over ds} - {2t\over m}
    {\xi_i\wedge\xi_i'\over(\xi_i',\xi_i')}.
$$

This characterization of semi-discrete isothermic surfaces
is an immediate consequence of \:def.sdi: in conjunction
with \:thm.gauge:.
To obtain the Darboux and Calapso transformations only the
flatness assertion is required, which already follows
from the simpler statement of \:thm.dtsgauge::
then the following characterizations of the Darboux and Calapso
transformations follow directly from their definitions and
the corresponding characterizations \:thm.darboux: and
\:def.calapso: for polarized curves,
cf \:thm.ttd:.

\proclaim\@thm.sditrafo:Thm 3.9@.
Let $\<\xi>:\Sigma\to S^n$ be semi-discrete isothermic with
isothermic loop of connections $(\Gamma^t,\nabla^t)_{t\in\R}$.
Then
\parindent=2em
\item{\rm(i)} the Darboux transforms $\<\hat\xi>$ of $\<\xi>$,
 with respect to a parameter $\mu$, are given by
 $(\Gamma^\mu,\nabla^\mu)$-parallel sections,
 that is, by solutions of Darboux's linear system
 $$
   \forall(ij)\in E:\Gamma^\mu_{ij}\hat\xi_j = \hat\xi_i
    \enspace{\sl and}\enspace
   \forall i\in V:\nabla^\mu_i\hat\xi_i = 0;
 $$
\item{\rm(ii)} the Calapso transforms $\<\xi^t>$ of $\<\xi>$
 are given as images of $\<\xi>$ under trivializing gauge
 transformations,
 $$
   \<\xi^t> = \<T^t\xi>,
    \enspace{\sl where}\enspace
   T^t\cdot(\Gamma^t,\nabla^t) = (\id,{d\over ds})
    \enspace\Leftrightarrow\enspace\cases{
   \forall(ij)\in E:
    T^t_i\circ\Gamma^t_{ij} = T^t_j, \cr
   \forall i\in V:
    T^t_i\circ\nabla^t_i = {d\over ds}\circ T^t_i. \cr}
 $$

At this point the full integrable theory of isothermic surfaces
is available and corresponding results in the semi-discrete
setting can be obtained in a completely analogous way as
in the smooth or fully discrete settings.
For example,
we may now define {\em special isothermic surfaces\/} as those
semi-discrete isothermic surfaces that admit a {\em polynomial
conserved quantity\/}, i.e., a polynomial map
$$
  p(t) = zt^d+yt^{d-1}+\dots+q:\Sigma\to\R^{n+1,1}
   \enspace{\rm with}\enspace
  \forall t\in\R:(\Gamma^t,\nabla^t)\,p(t) = 0,
$$
cf [\:busa12:, Defs 2.1 \& 2.3] and [\:bjrs15:, Defs 1 \& 2].
Wheeling out the conserved quantity conditions
$$
  (zt+q)_i=\Gamma_{ij}^t(zt+q)_j
   \enspace{\rm and}\enspace
  0 = \dts(zt+q)
$$
in the case of a {\em linear conserved quantity\/} $p(t)=zt+q$
yields $q\equiv const$, $z\perp\xi$ and, cf [\:bjrs15: (2.3)],
$$
  d_{ij}z = {1\over\mu_{ij}}(\pi_i-\pi_j)q
   \enspace{\rm and}\enspace
  z' = {2\over m(\xi',\xi')}\{(q,\xi)\xi'-(q,\xi')\xi\},
$$
where $\pi_i$ and $\pi_j$ denote the projections onto
$\<\xi_i>$ resp $\<\xi_j>$, and $d_{ij}z=z_j-z_i$ the
discrete derivative of $z$, as before;
note that
$
  \pi_j(d_{ij}z+{1\over\mu_{ij}}(\pi_j-\pi_i)q)
  = \pi_j(z_j+{1\over\mu_{ij}}\,q)
$
since $z\perp\xi$.
Thus normalizing $\xi=x$ so that $(x,q)\equiv -1$ we learn
that $z$ is, up to scale, a Christoffel transform of $x$
in $\R^{n+1,1}$:
$$
  z' = {2\over mx'}
   \enspace{\rm and}\enspace
  d_{ij}z = {2\over\mu_{ij}d_{ij}x}.
$$
As $(\Gamma^t,\nabla^t)$ are metric connections the (real)
polynomial $|zt+q|^2$ has constant coefficients,
in particular, we may without loss of generality assume
$|z|^2\equiv 1$ as long as $|z|\neq 0$, that is, $p(t)=zt+q$
to be a {\em normalized linear conserved quantity\/}.
Conversely, $p(t)=zt+q$ yields a (normalized) linear conserved
quantity of $\<\xi>=\<x>$ as soon as
$$
  z:\Sigma\to S^{3,1} := \{y\in\R^{n+1,1}\,|\,(y,y)=1\}
   \enspace{\rm with}\enspace
  H := -(z,q) = const
$$
is (up to scale) a Christoffel transform of
$$
  x:\Sigma\to Q^3
  = \{y\in{\cal L}^{n+1}\,|\,(y,q)=-1\}
  \subset\R^{n+1,1}.
$$

On the other hand, the two surfaces $x,z:\Sigma\to\R^{n+1,1}$
form a Christoffel pair if and only if their
{\em mixed area element\/}
$$
  A(x,z)ds := {1\over 2}
  (x'_{ij}\wedge d_{ij}z + z'_{ij}\wedge d_{ij}x)\,ds
  \equiv 0,
\eqno ma$$
with $z_{ij}={z_i+z_j\over 2}$, as before;
cf [\:muwa13:, Sect 5] and [\:mu15:, Sect 1.5],
see also [\:bosu08:, Thm 4.42] or [\:bjl14:, Lemma 2.3].

Namely, suppose that $z$ and $x$ define {\em parallel nets\/},
that is, with suitable functions $a_i$ and $\alpha_{ij}$
we have
$$
  z_i' = a_ix_i'
   \enspace{\rm and}\enspace
  d_{ij}z = \alpha_{ij}d_{ij}x.
$$
For regularity we assume $d_{ij}x\not\,\parallel x_i',x_j'$
and that $x$ and $z$ are not homothetic.
Integrability of $z$ then yields
$$
  (d_{ij}z)'=d_{ij}(z')
   \enspace\Leftrightarrow\enspace
  (\alpha_{ij}-a_i)\,x_i' - (\alpha_{ij}-a_j)\,x_j'
  = \alpha_{ij}'d_{ij}x.
\eqno intpn1$$
In particular, we learn that $x_i'$, $x_j'$ and $d_{ij}x$
must be linearly dependent,
that is, $x$ must be a {\em conjugate net\/} in $\R^{n+1,1}$,
cf [\:psbswbw08:, Sect 2];
hence adjacent curves $x_i$ and $x_j$ form Ribaucour pairs,
since $x$ maps into the conformal $n$-sphere,
and $x$ is {\em circular\/}
or a semi-discrete {\em curvature line net\/},
cf [\:muwa13:, Def 1.1].
Now, we use (\:intpn1:) to analyze vanishing of the mixed area,
$$
  A(x,z)_{ij} = 0
   \enspace\Leftrightarrow\enspace
  (\alpha_{ij}+a_i)\,x_i' + (\alpha_{ij}+a_j)\,x_j'
  = 0 \bmod d_{ij}x
   \enspace\Leftrightarrow\enspace
  \alpha_{ij}^2 = a_ia_j.
$$
We conclude that $a_ia_j>0$ and that there is a function
$\nu:\Sigma\to\R^\times$ with $a_i=\mp{1\over\nu_i^2}$
and $\alpha_{ij}={1\over\nu_i\nu_j}$.
Without loss of generality we may assume $a_i=-{1\over\nu_i^2}$
to recover that vanishing of the mixed area is a characterization
of the Koenigs duality of conjugate nets
from [\:muwa13:, Def 3.1],
up to a sign:
$$
  z_i' = -{1\over\nu_i^2}\,x_i'
   \enspace{\rm and}\enspace
  d_{ij}z = {1\over\nu_i\nu_j}\,d_{ij}x.
\eqno koenigs$$

However, for a semi-discrete curvature line net,
existence of a Koenigs dual implies isothermicity,
and the Koenigs and Christoffel dualities coincide.
To see this we use again integrability (\:intpn1:) of $z$,
which now reads
$$
  0 = ({1\over\nu_i\nu_j})'d_{ij}x
  + ({1\over\nu_i}+{1\over\nu_j})\,d_{ij}({x'\over\nu})
  = ({1\over\nu_i\nu_j})'d_{ij}x + {2\over\nu_i\nu_j}d_{ij}x'
  + ({1\over\nu_j}-{1\over\nu_i})\,({x'\over\nu})_{ij}.
\eqno intpn2$$
Since $x$ is isotropic, $(x,x)\equiv 0$,
we have $x_{ij}\perp d_{ij}x$,
hence using $d_{ij}({x'\over\nu})\parallel d_{ij}x$
from (\:intpn2:) we obtain
$$
  (d_{ij}x,({x'\over\nu})_{ij})
  = d_{ij}(x,{x'\over\nu})
  - (x_{ij},d_{ij}({x'\over\nu}))
  = d_{ij}(x,{x'\over\nu})
  = 0
   \enspace\Rightarrow\enspace
  d_{ij}({x'\over\nu}) \perp ({x'\over\nu})_{ij};
$$
further, taking inner product with $d_{ij}x$ in (\:intpn2:)
the right term yields the second of the equations
$$
  d_{ij}({1\over m}) = 0
   \enspace{\rm and}\enspace
  ({1\over\mu_{ij}})' = 0,
   \enspace{\rm where}\enspace
  {1\over m} := {(x',x')\over\nu^2}
   \enspace{\rm and}\enspace
  {1\over\mu_{ij}} := {2(x_i,x_j)\over\nu_i\nu_j},
$$
cf (\:sdicr:): now
(\:koenigs:) turns into Christoffel's equation (\:christoffel:)
and
(\:intpn2:) shows that $x_j$ is a $\mu_{ij}$-Darboux transform
of $x_i$, that is, isothermicity of $x$, as
$$
  0 = ({1\over\nu_i\nu_j})'d_{ij}x
  + ({1\over\nu_i}+{1\over\nu_j})\,d_{ij}({x'\over\nu})
  = ({1\over\nu_i}+{1\over\nu_j})\{
    {1\over\nu_j}\dts\big|_{t=\mu_{ij}}x_j
    - {1\over\nu_i}{(x_i',x_j)\over(x_i,x_j)}\,x_j
    \},
$$
where we use
${(x_i',x_j)\over\nu_i}={(x_i,x_j')\over\nu_j}$
from $({x'\over\nu})_{ij}\perp d_{ij}x$
and
$ {\nu_i'\over\nu_i}+{\nu_j'\over\nu_j}
  = {(x_i',x_j)+(x_i,x_j')\over(x_i,x_j)}
$
from the definition of $\mu_{ij}$.

Note that these results descend to Euclidean ambient geometry;
however, corresponding proofs in a purely Euclidean setting
will require some arguments to be adapted,
for example,
by using the tangent cross ratio instead of isotropy of $x$
to obtain $m$ and $\mu$ as functions of one parameter.

Using the vanishing of the mixed area (\:ma:) we then arrive
at characterizations of semi-discrete surfaces of constant
{\em mixed area mean curvature\/} in space forms,
$$
  H := -{A(x,n)\over A(x,x)},
   \enspace{\rm where}\enspace
  x:\Sigma\to Q^3
   \enspace{\rm and}\enspace
  n:\Sigma\to P^3 = \{y\in S^{3,1}\,|\,y\perp q\}
   \enspace{\rm with}\enspace
  n\perp x
\eqno cmc$$
defines a suitable {\em tangent plane congruence\/} for $x$
in $Q^3$,
cf [\:muwa13:, Prop 5.2] and [\:mu15:, Def 8] for the Euclidean
case and see [\:bjr14:, Def 2.3] for arbitrary ambient space
form geometries in the fully discrete case:
namely, consider
$$
  z = n + H\,x
   \enspace{\rm with}\enspace
  H = -(z,q)
$$
to see that $(x,n):\Sigma\to Q^3\times P^3$ has
 constant mean curvature
if and only if $x$ has a
 suitable Christoffel dual $z:\Sigma\to S^{3,1}$,
if and only if $\<\xi>=\<x>$ has a
 normalized linear conserved quantity $p(t)=zt+q$,
cf [\:bjl14:, Lemma~4.1] and [\:bjr14:, Thm~2.8]
resp [\:busa12:, Prop~2.5] and [\:bjrs08:, Sect~5].

\section References

\bgroup\frenchspacing\parindent=2em

\refitem bi04
 L Bianchi:
 {\it Il teorema di permutabilit\`a per le trasformazioni
  di Darboux delle superficie isoterme\/};
 Rend Acc Naz Lincei 13, 359--367 (1904)

\refitem bopi96
 A Bobenko, U Pinkall:
 {\it Discrete isothermic surfaces\/};
 J reine angew Math 475, 187--208 (1996)

\refitem boje01
 A Bobenko, U Hertrich-Jeromin:
 {\it Orthogonal nets and Clifford algebras\/};
 T\^ohoku Math Publ 20, 7--22 (2001)

\refitem bosu07a
 A Bobenko, Y Suris:
 {\it On organizing principles of discrete differential 
  geometry. Geometry of spheres\/};
 Russian Math Surveys 62, 1--43 (2007)

\refitem bosu08
 A Bobenko, Y Suris:
 {\it Discrete differential geometry.
  Integrable structure\/};
 Grad Stud Math 98, Amer Math Soc, Providence RI (2008)

\refitem bjl14
 A Bobenko, U Hertrich-Jeromin, I Lukyanenko:
 {\it Discrete constant mean curvature nets in space forms:
  Steiner's formula and Christoffel duality\/};
 Discr Comp Geom 52, 612--629 (2014)

\refitem bu06
 F Burstall:
 {\it Isothermic surfaces: conformal geometry,
  Clifford algebras and integrable systems\/};
 AMS Stud Adv Math 36, 1--82 (2006)

\refitem bdpp11
 F Burstall, N Donaldson, F Pedit, U Pinkall:
 {\it Isothermic submanifolds of symmetric $R$-spaces\/};
 J reine angew Math 660, 191--243 (2011)

\refitem bjrs08
 F Burstall, U Hertrich-Jeromin, W Rossman, S Santos:
 {\it Discrete surfaces of constant mean curvature\/};
 RIMS Kokyuroku 1880, 133--179 (2014)

\refitem busa12
 F Burstall, S Santos:
 {\it Special isothermic surfaces of type $d$\/};
 J London Math Soc 85, 571--591 (2012)

\refitem bjrs15
 F Burstall, U Hertrich-Jeromin, W Rossman, S Santos:
 {\it Discrete special isothermic surfaces\/};
 Geom Dedicata 174, 1--11 (2015)

\refitem bjr14
 F Burstall, U Hertrich-Jeromin, W Rossman:
 {\it Discrete linear Weingarten surfaces\/};
 EPrint arXiv: math.DG/1406.1293 (2014)

\refitem caiv98
 A Calini, T Ivey:
 {\it B\"acklund transformations and knots of
  constant torsion\/};
 J Knot Theor Ramifications 7, 719--746 (1998)

\refitem ci97
 J Cie\'sli\'nski:
 {\it The cross ratio and Clifford algebras\/};
 Adv Appl Clifford Alg 7, 133--139 (1997)

\refitem da99
 G Darboux:
 {\it Sur les surfaces isothermiques\/};
 Ann Sci \'Ec Norm Sup 16, 491--508 (1899)

\refitem jepe97
 U Hertrich-Jeromin, F Pedit:
 {\it Remarks on the Darboux transform of
  isothermic surfaces\/};
 Doc Math 2, 313--333 (1997)

\refitem imdg
 U Hertrich-Jeromin:
 {\it Introduction to M\"obius differential geometry\/};
 London Math Soc Lect Note Series 300,
  Cambridge Univ Press, Cambridge (2003)

\refitem ho09
 T Hoffmann:
 {\it Discrete differential geometry of curves and surfaces\/};
 COE Lect Notes 18, Kyushu Univ (2009)

\refitem muwa13
 C M\"uller, J Wallner:
 {\it Semi-discrete isothermic surfaces\/};
 Res Math 63, 1395--1407 (2013)

\refitem mu15
 C M\"uller:
 {\it Semi-discrete constant mean curvature surfaces\/};
 Math Z 279, 459--478 (2015)

\refitem psbswbw08
 H Pottmann, A Schiftner, P Bo, H Schmiedhofer, W Wang,
  N Baldassini, J Wallner:
 {\it Freeform surfaces from single curved panels\/};
 ACM Trans Graph 27, 76 (2008)

\refitem roya12
 W Rossman, M Yasumoto:
 {\it Weierstrass representation for semi-discrete minimal
  surfaces, and comparison of various discretized catenoids\/};
 J Math Industry 4, 109--118 (2012)

\egroup


\vskip3em\vfill
\bgroup\fn[cmr7]\baselineskip=8pt
\def\addwd{\hsize=.36\hsize}
\def\udo{\vtop{\addwd
 U Hertrich-Jeromin\\
 Vienna University of Technology\\
 Wiedner Hauptstra\ss{}e 8--10/104\\
 A-1040 Vienna (Austria)\\
 Email: udo.hertrich-jeromin@tuwien.ac.at
 }}
\def\fran{\vtop{\addwd
 F Burstall\\
 Department of Mathematical Sciences\\
 University of Bath\\
 Bath, BA2 7AY (United Kingdom)\\
 Email: f.e.burstall@bath.ac.uk
 }}
\def\wayne{\vtop{\addwd
 W Rossman\\
 Department of Mathematics\\
 Kobe University\\
 Rokko, Kobe 657-8501 (Japan)\\
 Email: wayne@math.kobe-u.ac.jp
 }}
\def\christian{\vtop{\addwd
 C M\"uller\\
 Vienna University of Technology\\
 Wiedner Hauptstra\ss{}e 8--10/104\\
 A-1040 Vienna (Austria)\\
 Email: cmueller@geometrie.tuwien.ac.at
 }}
\hbox to \hsize{\hfil \fran \hfil \udo \hfil}\vskip 3ex
\hbox to \hsize{\hfil \christian \hfil \wayne \hfil}
\egroup
\bye